\documentclass[11pt]{article}
\usepackage{geometry}
\usepackage{setspace}
\usepackage[OT1]{fontenc}
\usepackage{amssymb,amsmath,amsthm}

\newcommand\X{\mathbf X}
\newcommand\x{\mathbf x}
\def\Im{{\rm\,Im} }
\newcommand\bR{{\mathbb R}}

\newcommand\e{{\mathbb E}}
\newcommand\p{{\mathbb P}}
\newcommand\var{{\rm Var}}
\newcommand\tr{{\mathrm{tr}}}

\newcommand\bC{{\mathbb C}}

\newcommand\pto{\stackrel p\rightarrow}

\newtheorem{lemma}{Lemma}[section]%
\newtheorem{theorem}{Theorem}[section]%
\newtheorem{proposition}{Proposition}[section]%

\begin{document}
\begin{center}
\Large A short proof of the Marchenko-Pastur theorem.  
\end{center}
\begin{center}
\large Pavel~Yaskov\footnote{Steklov Mathematical Institute, Russia\\
Moscow Institute of Physics and Technology\\
 e-mail: yaskov@mi.ras.ru\\Supported
    by RNF grant 14-21-00162 from the Russian Scientific Fund.}
 \end{center}

\begin{abstract}
We prove the Marchenko-Pastur theorem for random matrices with i.i.d. rows and a general dependence structure  within the rows by  a simple modification of the standard Cauchy-Stieltjes resolvent method. 


\begin{center}
\end{center}
\end{abstract}

\section{Introduction}
Let $\X_{pn}$ be a $p\times n$ random matrix whose columns $\{\x_{pk}\}_{k=1}^n$ are i.i.d. copies of some random vector  $\x_p$ in $\bR^p$ for all $p,n\geqslant1$.  All random elements are defined on the same probability space. The object of our study is $\mu_{pn}$, the empirical spectral distribution (ESD) of $n^{-1}\X_{pn}\X_{pn}^\top$. Here ESD of a $p\times p$ real symmetric matrix $A$ is defined by 
\[\mu=\frac{1}{p}\sum_{i=1}^p \delta_{\lambda_{i}},\]
where $\delta_\lambda$ stands for the Dirac mass at $\lambda\in\bR$ and $\lambda_1\leqslant\ldots\leqslant\lambda_p$ are eigenvalues of $A.$ 

Recall that the Marchenko-Pastur law $\mu_c$ with parameter $c>0$ is the probability distribution
\[(1-1/c)^+\delta_0+\frac{\sqrt{(b-x)(x-a)}}{2\pi c x}I(x\in [a,b])\,dx,\]
  where  $x^+=\max\{x,0\}$ for $x\in\bR,$  $a=(1-\sqrt{c})^2$, and  $b=(1+\sqrt{c})^2.$

The Marchenko-Pastur theorem states that, for any $p=p(n)$ with $p/n\to c>0$ as $n\to\infty$, \begin{equation}
\label{mp} \p(\text{$\mu_{pn}\to\mu_c$ in distribution, } n\to\infty)=1\end{equation} if each $\x_p$ has centred orthonormal entries $\{X_{pk}\}_{k=1}^p$  satisfying certain conditions.  

The standard conditions include the independence of $\{X_{pk}\}_{k=1}^p$ and the Lindeberg condition
\begin{equation}\label{l}
\lim_{p\to\infty}\frac{1}{p}\sum_{k=1}^p \e X_{pk}^2I( |X_{pk}|>\varepsilon \sqrt{p})=0\quad\text{for all }\varepsilon>0
\end{equation}
(see Theorem 3.10 in \cite{BS}). Bai and Zhou \cite{BZ}, Pastur and Pajor \cite{PP}, and Pastur and Scherbina \cite{PS} (see Theorem 19.1.8) proved the Marchenko-Pastur theorem, assuming that $\var(\x_p^\top A_p\x_p/p)\to 0$, $p\to\infty,$ for all sequences of  $p\times p$ complex  matrices $A_p$ with uniformly bounded spectral norms $\|A_p\|$. If entries of $\x_p$ are independent, this assumption is much stronger than \eqref{l}. 

In this note we give a short proof of the  Marchenko-Pastur theorem  under weaker conditions that cover all mentioned results.

\section{Main results}
Consider the following assumption.\\
(A) $(\x_p^\top A_p\x_p-\tr(A_p))/p\pto 0$ as $p\to\infty$ for all sequences of  $p\times p$ complex matrices $A_p$ with uniformly bounded spectral norms $\|A_p\|$.
\begin{theorem}\label{t1}
If $(A)$ holds, then \eqref{mp} holds.
\end{theorem} 

If entries of $\x_p$ are orthonormal, then $\e \x_p^\top A_p\x_p=\tr(A_p)$ and the assumption considered in \cite{BZ}, \cite{PP}, \cite{PS} (see Introduction) is stronger than (A).  In addition, we have the following proposition.

\begin{proposition}\label{p}
Let $\{X_{pk}\}_{k=1}^p$ be independent random variables with $\e X_{pk}=0,$ $\e X_{pk}^2=1$ for each $p\geqslant 1$. Then  \eqref{l} holds if and only if $(A)$ holds for $\x_p=(X_{p1},\ldots,X_{pp})$, $p\geqslant 1$.
\end{proposition} 

 Assumption (A) also covers the case, where entries of $\x_p$ are orthonormal infinite linear combinations (in $L_2$) of some i.i.d. random variables $\{\varepsilon_k\}_{k=1}^\infty$ with $\e\varepsilon_k=0$ and $\e\varepsilon_k^2=1$ (see Corollary 4.9 in arXiv:1410.5190).
  
\noindent{\bf Remark.} We get an equivalent reformulation of (A) if we replace complex matrices by real symmetric positive semi-definite matrices (see Lemma 5.3 in arXiv:1410.5190).

\section{Proofs}
\begin{proof}[\bf Proof of Theorem \ref{t1}]
 We will use the Cauchy-Stieltjes transform method. By the Stieltjes continuity theorem (e.g., see Exercise 2.4.10(i) in \cite{T}), we only need to show that $s_n(z)\to s(z)$ a.s. for all $z\in\bC$ with $\Im(z)>0$, where $s_n=s_n(z)$ and $s=s(z)$ are
the Stieltjes transforms of $\mu_{pn}$ and $\mu_c$ defined by
\[s_n(z)=\int_\bR \frac{\mu_{pn}(d\lambda)}{\lambda-z}\quad\text{and}\quad s(z)=\int_\bR \frac{\mu_c(d\lambda)}{\lambda-z}.\]
By the definition of $\mu_{pn}$, $s_n(z)=\tr(n^{-1}\X_{pn}\X_{pn}^\top-zI_p)^{-1}/p$ for the $p\times p$ identity matrix $I_p$.

Fix any $z\in\bC$ with $v=\Im(z)>0$.  By the standard martingale argument (e.g., see Step 1 in the proof of Theorem 1.1 in \cite{BZ} or Lemma 4.1 in \cite{A}), we derive that $s_n(z)-\e s_n(z)\to 0$ a.s.  We finish the proof by checking that $\e s_n(z)\to s(z).$ We need a technical lemma.
\begin{lemma} \label{l1} Let  $C$ be a real symmetric positive semi-definite $p\times p$ matrix and $x\in\mathbb \bR^p$. If  $z\in\mathbb C$ is such that $v=\Im(z)>0,$ then $(1)$ $\|(C-zI_p)^{-1}\|\leqslant 1/v$, $(2)$ 
$|\tr(C+xx^\top-zI_p)^{-1}-\tr(C-zI_p)^{-1}|\leqslant 1/v,$ $(3)$ $|x^\top (C+xx^\top-z I_p)^{-1}x)|\leqslant 1+|z|/v$,  $(4)$ $\Im(z+z\tr(C-z I_p)^{-1})\geqslant v $ and $\Im(\tr(C-z I_p)^{-1})>0$, $(5)$  $\Im(z+zx^\top(C-z I_p)^{-1}x)\geqslant v$.
\end{lemma}
All bounds in Lemma \ref{l1} are well-known. Part (1) can be proved by diagonalizing $C.$ Part (2) is given in Lemma 2.6 in \cite{BS1}. Part (3) follows from the Sherman-Morrison formula and Part (5), since
  \[x^\top (C+xx^\top-z I_p)^{-1}x=x^\top (C-z I_p)^{-1}x-\frac{(x^\top (C-z I_p)^{-1}x)^2}{1+x^\top (C-z I_p)^{-1}x}=1-\frac{z}{z+zx^\top (C-z I_p)^{-1}x}.\]
  Parts (4)--(5) can be checked by showing that $\Im(\tr((1/z)C-I_p)^{-1})\geqslant 0 $ and  $\Im(x^\top((1/z)C- I_p)^{-1}x)\geqslant 0$.

Take $\x_p=\x_{p,n+1}$ to be independent of the matrix $\X_{pn}$ and distributed as its columns $\{\x_{pk}\}_{k=1}^n$. Define \[A_n=\X_{pn}\X_{pn}^\top=\sum_{k=1}^n\x_{pk}\x_{pk}^\top\quad \text{and }\quad B_n=A_n+\x_{p}\x_{p}^\top=\sum_{k=1}^{n+1}\x_{pk}\x_{pk}^\top.\]
By Lemma \ref{l1}(1), $B_n -znI_p$ is non-degenerate and
\[p=\tr\big((B_n -znI_p)(B_n-znI_p)^{-1}\big)=\sum_{k=1}^{n+1}\x_{pk}^\top(B_n -znI_p)^{-1}\x_{pk}-zn\,\tr(B_n-znI_p)^{-1}.\]
Taking expectations and using the exchangeability of $\{\x_{pk}\}_{k=1}^{n+1}$,
\begin{align}\label{e1}
p=&(n+1)\e \x_p^\top(B_n -znI_p)^{-1}\x_p-zn\,\e \tr(B_n-znI_p)^{-1}.
\end{align}

Define $S_n(z)=\tr(A_n -znI_p)^{-1}$ and note that $S_n(z)=(p/n)s_n(z)$.   By Lemma \ref{l1}(2)-(3), 
\[\e\tr(B_n-znI_p)^{-1}=\e S_n(z)+O(1/n)\quad\text{and}\quad \e \x_p^\top(B_n -znI_p)^{-1}\x_p=O(1).\]
Moreover,  we will show below that
\begin{equation}\label{m}
\e \x_p^\top(B_n -znI_p)^{-1}\x_p=\frac{\e S_n(z)}{1+\e S_n(z)}+o(1)
\end{equation}

Suppose for a moment that \eqref{m} holds (and $p/n=c+o(1)$). Then \eqref{e1} reduces to
\[\frac{\e S_n(z)}{1+\e S_n(z)}-z\e S_n(z)=c+o(1).\]
By  (1) and (4) in Lemma \ref{l1}, $(\e S_n(z))_{n=1}^\infty$  is a bounded sequence with $\Im(\e S_n(z))>0$. It is straightforward to check that the limiting quadratic equation  $S/(1+S)-zS=c$ or $zS^2+(z-1+c)S+c=0$ has a unique solution $S=S(z)$ with $\Im(S(z))\geqslant 0$ when $\Im(z)>0.$  As a result, any converging subsequence of bounded sequence $(\e S_n(z))_{n=1}^\infty$ tend to $S(z)$. This implies that $\e S_n(z)=(p/n)\e s_n(z)\to S(z).$

One can also show that $S(z)=cs(z)$ is the above unique solution, where $s(z)$ is the Stieltjes transform of the Marchenko-Pastur law  (see Remark 1.1 in \cite{BZ}). Combining all above relations, we conclude that $s_n(z)\to s(z)$ a.s. 

To finish the proof, we only need to check \eqref{m}.  By the Sherman-Morrison formula,
\[\x_p^\top(B_n -znI_p)^{-1}\x_p=\x_p^\top(A_n+\x_p\x_p^\top -znI_p)^{-1}\x_p=\frac{\x_p^\top(A_n -znI_p)^{-1}\x_p}{1+\x_p^\top(A_n -znI_p)^{-1}\x_p}.\]
Using Lemma \ref{l1}(1), $(A)$, and the independence of $\x_p$ and $A_n$, we get $\x_p^\top(A_n-znI_p)^{-1}\x_p-S_n(z)\pto 0.$ In addition, as it is shown above, \[S_n(z)-\e S_n(z)=(p/n)(s_n(z)-\e s_n(z))\pto0.\] 
Hence, Lemma \ref{l1}(4)-(5) and  inequality $|1+w|\geqslant \Im(z+zw)/|z|,$  $w\in\bC$, yield
\[\bigg|\frac{\x_p^\top(A_n -znI_p)^{-1}\x_p}{1+\x_p^\top(A_n -znI_p)^{-1}\x_p}-\frac{\e S_n(z)}{1+\e S_n(z)}\bigg|\leqslant 
\frac{|z|^2}{v^2}|\x_p^\top(A_n -znI_p)^{-1}\x_p-\e S_n(z)|\pto 0.\]
Finally, \eqref{m} follows from Lebesgue's dominated convergence theorem and  Lemma \ref{l1}(3).
\end{proof}
\begin{proof}[\bf Proof of Proposition \ref{p}.] For each $p\geqslant 1,$ let $A_p=\big(a_{kj}^{(p)}\big)_{k,j=1}^p$ be a complex $p\times p$ matrix  with $\|A_p\|\leqslant 1$. If $D_p$ is a diagonal matrix with diagonal entries $\big(a_{kk}^{(p)}\big)_{k=1}^p$, then
\[\e\Big|\x_p^\top (A_p-D_p)\x_p\Big|^2\leqslant 
2\e\Big|\sum_{1\leqslant k<j\leqslant p}a_{kj}^{(p)}X_{pk}X_{pj}\Big|^2+2\e\Big|\sum_{1\leqslant j<k\leqslant p}a_{kj}^{(p)}X_{pk}X_{pj}\Big|^2=2\sum_{j\neq k} |a_{jk}^{(p)}|^2\leqslant 4 \tr(A_pA_p^*),
\]
where $A_p^*$ is the complex conjugate of $A_p.$ By the definition of the spectral norm, $\tr(A_pA_p^*)\leqslant \|A_p\|^2 p.$ Thus,
\[\frac{\x_p^\top (A_p-D_p)\x_p}{p}\pto 0.\]

To finish the proof, we   need to show that \eqref{l} holds if and only if 
\begin{equation}\label{e}
\frac{1}{p}\sum_{k=1}^p a_k^{(p)}(X_{pk}^2-1)\pto0\text{\quad for any triangular array $\{a_k^{(p)},1\leqslant k\leqslant p,p\geqslant 1\}$ with $|a_k^{(p)}|\leqslant 1$.}
\end{equation}

Let \eqref{e}  hold. Then $Z_p=p^{-1}\sum_{k=1}^p X_{pk}^2\pto 1$. Since $\e Z_p=1$ and $Z_p\geqslant 0$ a.s., the Vitali convergence theorem  yields $\e|Z_p-1|\to 0$. Using inequalities $p^{-1}\sum_{k=1}^p \e |X_{pk}^2-1|\leqslant 2,$ $p\geqslant 1$, we derive from \cite{H} that 
\[\frac{1}{p}\sum_{k=1}^p \e|X_{pk}^2-1|I(|X_{pk}^2-1|>\varepsilon p)\to0 \quad\text{for all $\varepsilon>0$}.\]
Obviously, this is equivalent to \eqref{l}.

Let \eqref{l} hold. By the Marcinkiewicz--Zygmund inequality, there is a universal constant $C>0$ such that
\[\e\Big|\frac{1}{p}\sum_{k=1}^p a_k^{(p)}(X_{pk}^2-1)\Big|\leqslant \frac{C}{p}\e\Big(\sum_{k=1}^p (X_{pk}^2-1)^2\Big)^{1/2}.\]
Using \eqref{l}, Jensen's inequality, and $\sqrt{x+y}\leqslant \sqrt{x}+\sqrt{y},$ $x,y\geqslant 0$, we get
\begin{align*}
\e\Big(\sum_{k=1}^p (X_{pk}^2-1)^2\Big)^{1/2}\leqslant&
\Big(\sum_{k=1}^p \e(X_{pk}^2-1)^2I(|X_{pk}^2-1|\leqslant\varepsilon p)\Big)^{1/2}+
\sum_{k=1}^p \e|X_{pk}^2-1|I(|X_{pk}^2-1|>\varepsilon p)\\
\leqslant& p\sqrt{2\varepsilon} +o(p)
\end{align*}
for all $\varepsilon>0,$ where we also applied the bound $\e(X_{pk}^2-1)^2I(|X_{pk}^2-1|\leqslant\varepsilon p)\leqslant \varepsilon p\e|X_{pk}^2-1|\leqslant 2 \varepsilon p$. Therefore, 
\[\varlimsup_{p\to\infty}\e\Big|\frac{1}{p}\sum_{k=1}^p a_k^{(p)}(X_{pk}^2-1)\Big|\leqslant \sqrt{2\varepsilon}.\]
Tending $\varepsilon$ to zero, we get \eqref{e}.\end{proof}

\end{document}